\documentclass[a4paper, reqno]{amsart}

\title[{}]
{Nesting behind $\hat{Z}$-invariants}

\author{Shoma Sugimoto}
\address{Division of Physics,
Mathematics and Astronomy, California Institute of Technology}
\email{shomasugimoto361@gmail.com}

\usepackage{latexsym}
\usepackage{amsmath}
\usepackage{amsfonts}
\usepackage{amssymb}
\usepackage{amsxtra}
\usepackage{mathrsfs}
\usepackage{eucal}
\usepackage[all]{xy}
\usepackage{color}
\usepackage{braket}
\usepackage{graphics, setspace}
\usepackage{etoolbox, textcomp}
\usepackage{comment}
\usepackage{changepage}
\usepackage{xcolor}
\definecolor{rouge}{rgb}{0.85,0.1,.4}
\definecolor{bleu}{rgb}{0.1,0.2,0.9}
\definecolor{violet}{rgb}{0.7,0,0.8}
\usepackage[colorlinks=true,linkcolor=bleu,urlcolor=violet,citecolor=rouge]{hyperref}
\usepackage{cleveref}
\usepackage{lscape}
\usepackage{tikz}
\usepackage{amsthm}
\usetikzlibrary{graphs}
\usetikzlibrary {arrows.meta}
\usepackage{enumitem}

\theoremstyle{definition}

\newtheorem*{fundamental-problem}{\underline{Fundamental Problem (FP)}}

\newtheorem*{guiding-principle}{\underline{Fundamental Problem from the new perspective}}

\numberwithin{equation}{section}

\newcommand{\Z}{\mathbb{Z}}

\newcommand{\id}{\operatorname{id}}

\newcommand{\g}{\mathfrak{g}}

\newcommand{\sll}{\mathfrak{sl}}

\newcommand{\SL}{SL_2}

\newcommand{\Bits}[1]{\Z_2^{#1}}

\newcommand{\rt}[1]{\operatorname{rt}({#1})}

\begin{document}
\maketitle

\begin{abstract}
In the spirit of \cite{S4}, 
we propose an abelian categorification of $\hat{Z}$-invariants \cite{GPPV} for negative definite plumbed 3-manifolds.
It provides a blueprint for the expected dictionary \cite{CFGHP} between these $3$-manifolds and log VOAs; that is, the contribution from 3d $\mathcal{N}=2$ theory via 3d-3d correspondence
is encoded as recursive and binary deviations from semisimplicity in the abelian category of modules over the hypothetical log VOA, and is decoded by the 
recursive application of the theory of Feigin--Tipunin construction \cite{LS}.
In particular, the nested Weyl-type character formulas provide virtual generalized characters reconstructing the $\hat{Z}$-invariants.
\end{abstract}

\section{Introduction}
\label{introduction}
The $\hat{Z}$-invariants are quantum invariants of negative definite plumbed 3-manifolds introduced in \cite{GPPV}, not only $q$-expansions of the WRT invariants \cite{M2} but also rich examples of spoiled modular forms \cite{CCFGH}. 
The last fact with \cite{Zhu} predict the
existence of irrational vertex operator algebras (log VOAs) with $\hat{Z}$-invariants as their virtual generalized characters \cite{CFGHP}.
However, no examples of such log VOAs have been found so far except for the two easiest cases, 3- or 4-leg star graphs, and a new perspective is needed to construct and study the general cases.

In this paper, by refining and generalizing \cite{S4}, for any negative definite plumbed 3-manifolds, we propose a hypothetical but unified construction/research method for the corresponding log VOAs. 
Let us give a digest of the background (for more detail, see \cite{S4}). 
Recall that $\hat{Z}$-invariants were defined as follows \cite[(2.3)]{CCFGH}:
\begin{align}
\label{Zhat-invariant}
\hat{Z}=\int\frac{dx}{2\pi ix}F_{3d}(x)\Theta_{2d}(q,x),
\end{align}
where $F_{3d}(x)$ corresponds the contribution from 3d $\mathcal{N}=2$ theory via 3d-3d correspondence (note that it depends only on the base tree) and $\Theta_{2d}(q,x)$ the 2d boundary degrees of freedom. 
Mathematically, $\Theta_{2d}(q)$ is the positive definite lattice theta function associated with the inverse of the linking matrix of the plumbed graph, i.e. a modular form, while the modularity of $\hat{Z}$ is spoiled by the contribution of $F_{3d}(x)$. 
Hence it is expected that if we can construct a categorification of $F_{3d}(x)$ and place the lattice VOA in the category, it will provide the
free field realization of the desired log VOA.
As such, we consider the nested Feigin--Tipunin constructions \cite{S4}. 

The FT construction is a geometric construction of log VOAs corresponding to 3-leg star graphs above, which was introduced and conjectured in \cite{FT}. 
The author proved this conjecture \cite{S1} and have developed the geometric representation theory and its applications to the log VOAs \cite{CNS,LS,S1,S2}. 
Actually, the FT construction itself does not require any VOA structures and can be developed in a purely Lie algebraic setting \cite{LS}. 
It provides a clear perspective for the study of log VOAs; that is, by decomposing the study into the general theory of FT construction and its application to easier (e.g. lattice) VOAs, we can avoid the difficult log VOA-structures that have been a bottleneck.
In fact, the above works predict that non-VOA-specific/abelian categorical aspects of the representation theory of the corresponding log VOAs are essentially reduced to the theory of FT construction. 

On the other hand, in the case of star graphs and $\g=\sll_2$, by recursively applying the Weyl-type character formula of the FT construction \cite{LS} to the lattice theta functions above, the $\hat{Z}$-invariants are obtained \cite{S3,S4}.
Combining it with the VOA-independence of the FT above, we can see the following: 
for each negative definite plumbed graph, or rather the base tree, if we can develop in advance the theory of ``nested FT constructions" 
on an abstract abelian category,
it will not only provide the categorification of $F_{3d}(x)$ above (or a blueprint for the abelian category of modules over the hypothetical log VOA/W-algebra) but also a powerful reductive approach to log VOAs.
In \cite{S4}, the author proposed such an abelian categorification for all star graphs (unless otherwise noted, we assume that $\g=\sll_2$ from now on) and thus a glimpse into the expected log VOAs/W-algebras.
We not only extends \cite{S4} to all trees, but also explains how the 
``virtual generalized characters"  reconstruct the $\hat{Z}$. 
Because of the emphasis on motive and reasoning process, this paper does not adopt the ``Bourbaki style".

\noindent\textbf{Acknowledgements:}
\noindent
The author is supported by Caltech-Tsinghua joint postdoc Fellowships, and was supported by Beijing Natural Science Foundation grant IS24011 until June 30th 2025.
This research was initiated and completed while the author was a postdoc at YMSC, Tsinghua University.

\section{Preliminaries}
\subsection{Notations}
A \textit{tree} is an undirected connected acyclic graph.
In this paper, we only consider finite trees.
Let $T=(V,E)$ be a tree with the vertex set $V$ and edge set $E$.
A \textit{path} (with length $n$) of $T$ is a sequence of edges $(e_i=\{v_i,v_{i+1}\})_{1\leq i\leq n}$.
For $v\in V$, the \textit{degree} $\deg(v)$ of $v$ is $|\{e\in E\mid v\in e\}|$.
For $d\in\Z_{\geq 1}$, 
set 
\begin{align*}
V_{\geq d}:=\{v\in V\mid \deg(v)\geq d\},
\ \ 
V_{\leq d}:=\{v\in V\mid \deg(v)\leq d\}, 
\ \ 
V_d:=V_{\geq d}\cap V_{\leq d}.  
\end{align*}
The vertices in $V_1$ and $V_{\geq 3}$ are called \textit{leaves} and \textit{nodes} of $T$, respectively.
On the other hand, there is ambiguity about which vertex is considered the \textit{root} of a tree.
A tree with a specified root vertex is called a \textit{rooted tree}.
Let $T$ be a rooted tree with root $\rt{T}$.
For any $v\in V$, $v\not=\rt{T}$, there exists a unique path $p_v$ from $v$ to $\rt{T}$.
The \textit{parent} of a vertex $v\not=\rt{T}$ is the unique vertex adjacent to $v$ on the path $p_v$ above.
A \textit{child} of $v\in V$ is a vertex with the parent $v$.
An \textit{ascendant} of $v\not=\rt{T}$ is a vertex $v'\in p_v$.
A \textit{descendant} of $v$ is a vertex $v'$ such that $v\in p_{v'}$.
A \textit{sibling} of $v$ is a vertex with the same parent as $v$.
The \textit{depth} of $v$ is the length of $p_v$.
The \textit{height} of $v$ (resp. $T$) is the length of a longest path from $v$ (resp. $\rt{T}$) to a leaf.
It is well known that a tree has only one or two adjacent \textit{centers}, i.e. the vertices with minimal height.
We call $T$ a \textit{centered tree} if $\mathrm{rt}(T)$ is a center.

A \textit{plumbed graph} $\Gamma=(T,(w_v)_{v\in V})$ is the \textit{base tree} $T=(V,E)$ with integer weights $w_v\in\Z$ for each $v\in V$.
By applying the Dehn surgery to the weighted link associated with $\Gamma$, the \textit{plumbed $3$-manifold} $M(\Gamma)$ is defined.
We consider the case where the linking matrix $W$ of $\Gamma$ is negative definite. 
Then (a family of) $q$-series called the \textit{$\hat{Z}$-invariants} of $M(\Gamma)$ are defined \cite{GPPV}.
Various formulas for $\hat{Z}$-invariants of $M(\Gamma)$ have been given and studied, but the important for us is the following ``\textit{bosonic form}" \cite{CFGHP} given by \cite[Proposition 4.2]{M1}:
\begin{align}
\label{RHS of M1, p.10}
\sum_{e\in\{\pm 1\}^{V_{\geq 3}}}
(\prod_{v\in V_{\geq 3}}e_v^{\deg(v)-|\bar{v}|})
\sum_{\epsilon\in\{\pm 1\}^{V_1}}
(\prod_{i\in\bar{v}}\epsilon_i)
\sum_{n\in\Z_{\geq 0}^{V_{\geq 3}}}
q^{Q(e(n+\frac{\deg(v)-2}{2}+\sum_{i\in\bar{v}}\frac{\epsilon_i}{2w_i}))}
\prod_{v\in V_{\geq 3}}\binom{n_v+\deg(v)-3}{n_v},
\end{align}
where 
$\bar{v}=\{i\ |\ (i,v)\in E \}\cap V_1$, 
$Q(x):={}^txSx$ 
is the positive definite quadratic form associated with the $V_{\geq 2}\times V_{\geq 2}$ submatrix $S$ of $-W^{-1}$ and 
$ex:=((x_v)_{v\in V_{2}},(e_vx_v)_{v\in V_{\geq 3}})$.
The $Q(x)$ defines the positive definite lattice theta function $\Theta_{2d}(q)$ in \eqref{Zhat-invariant}.
Our goal is to recursively reconstruct \eqref{RHS of M1, p.10} in an abelian categorical setting.

The abelian categorical terminology is the same as \cite[Section 3.1]{S4}, but we quote it here for convenience.
Let $M$ be an object of an abelian category.
We call $M$ \textit{semisimple} if it is a direct sum of simple subobjects.
A \textit{Loewy series} of $M$ is a minimal length series among strictly ascending series of subobjects 
$0=M_0\subset M_1\subset\cdots\subset M_{l(M)}=M$
such that each successive quotient $M_k/M_{k-1}$ is semisimple. 
In this paper, the minimal length $l(M)$ is always finite.
We use the letters $\iota_k$ and $\pi_k$ for the embedding 
$\iota_k\colon M_k\hookrightarrow M$
and the projection 
$\pi_k\colon M\twoheadrightarrow M/M_k$, respectively.
Simple subobjects in each $M_k/M_{k-1}$ are called \textit{composition factors} (CFs), and by Jordan-H\"{o}lder theorem, they are unique up to isomorphisms.
For a fixed Loewy series of $M$, the \textit{Loewy diagram} is the diagram with $l(M)$ layers such that the $k$-th layer consists of the CFs of $M_k/M_{k-1}$.
In addition, if the embedding of each CF $X_k$ into $M_k/M_{k-1}$ is also fixed, 
the \textit{annotated Loewy diagram} (ALD) is obtained by adding arrows to the Loewy diagram according to the following rule: 
for each $1\leq k\leq l(M)$ and pair of CFs $(X_k,X_{k-1})$, $X_k\subset M_k/M_{k-1}$ and $X_{k-1}\subset M_{k-1}/M_{k-2}$,
add the arrow $X_k\rightarrow X_{k-1}$ if there exists a subquotient $X\subseteq M_k/M_{k-2}$ with a non-split short exact sequence $0\rightarrow X_{k-1}\rightarrow X\rightarrow X_k\rightarrow 0$.
In general, ALDs are not uniquely determined from mere Loewy series.
However, if all CFs are non-isomorphic each other, then the ALD is uniquely determined.
Let $C$ be the set of isomorphism classes of CFs in $M$.
By coloring each CF $X_k$ with its isomorphism class, we obtain the $C$-colored \textit{directed acyclic graph} (DAG) from an ALD.
Conversely, we can also consider the case where for a given $C$-colored DAG $Q$, there exists an object $M$ with an ALD represented by $Q$ in the above sense.
When an ALD of $M$ and a $C$-colored DAG $Q$ are identified in the above-mentioned sense, 
the structures of $M$ and $Q$ corresponds as follows:
a CF $X_k$ in the ALD of $M$ and the corresponding node in $Q$,
the subobject $\pi_k^{-1}(X_k)$ and the  induced subgraph of $Q$ generated by the node $X_k$,
the inclusion relation among $\pi_k^{-1}(X_k)$ and reachability in $Q$,
and
the isomorphism class of $X_k$ and the corresponding color in $C$,
etc. 
As mentioned above, if $Q$ is \textit{$C$-labeled} (i.e. the coloring is injective), then the ALD of $M$ is $Q$.
But even if not, we sometimes identify $M$ with its ALD/DAG if no confusion can arise.
The ALD/DAG $Q^\ast$ is obtained br reversing all arrows in $Q$.
The Grothendieck group of an abelian category $\mathcal{C}$ is denoted by $[\mathcal{C}]$.

\subsection{Review of \cite{S4}}
\label{quick review of S4}
Let $C_m:=\Z_2^m\times\Z_{\geq 0}$.
The elements in $\Z_2^m$ and $\Z_{\geq 0}$ are called \textit{bits} and \textit{depth}, respectively.
We regard $\Z_2=\{+,-\}$ as the Weyl group $(\{\id,\sigma_1\},\ast)$ of $\sll_2$.
For $\lambda=\lambda_1\cdots\lambda_m\in\Z_2^m$, set $|\lambda|:=|\{i\ |\ \lambda_i=-\}|$.
Also, for the partial order $\lambda\leq\mu\Leftrightarrow(\lambda_i=+\Rightarrow\mu_i=+)$ on $\Z_2^m$, set $|\mu-\lambda|:=|\{i\ |\ \lambda_i=-\ \land\ \mu_i=+\}|\in\Z_{\geq 0}$.
For $\lambda=\lambda_1\cdots\lambda_m\in\{+,-,\pm\}^m$, we identify $\lambda$ with the subset $\{\mu\in\Z_2^m\ |\ \lambda_i\in\Z_2\Rightarrow\mu_i=\lambda_i\}\subseteq\Z_2^m$.

Let us introduce two basic objects.
First, the $C_m$-colored $2m$-cube DAG $(\pm^m|\pm^m)$ is defined by 
\begin{align*}
&\operatorname{Node}(\pm^m|\pm^m)
:=\{
(\lambda|\mu)
\ |\ 
(\lambda,\mu)\in \Bits{m}\times\Bits{m}=\Bits{2m}
\},
\ \ 
b(\lambda|\mu)
:=\lambda\ast\mu\in\Z_2^m,
\ \ 
d(\lambda|\mu)
:=|\lambda\mu|\in\Z_{\geq 0},
\\
&\operatorname{Edge}(\pm^m|\pm^m)
:=
\{
(\lambda|\mu)\rightarrow(\lambda'|\mu')
\ |\ 
(|\lambda'-\lambda|=1
\ \land\ 
\mu=\mu')
\lor
(|\mu-\mu'|=1
\ \land\ 
\lambda=\lambda')
\}
\end{align*}
\cite[Definition 2.1]{S4}.
(Let's draw some pictures when $m$ is small.)
More generally, for subsets $D,E\subseteq\Z_2^m$, denote $(D|E)$ the corresponding induced subgraph.
In particular, for $\lambda\in\Z_2^m$ and $\bar\lambda:=-^m\ast\lambda$, the $m$-cube DAGs $(\lambda|\pm^m)$ and $(\pm^m|\bar\lambda)$ are considered.
For a $C_a$-colored DAG $Q$, we write the Cartesian product $Q\square(D|E)$ as $Q(D|E)$ if no confusion can arise.
For such products, the composition rules reflecting the hypercubic structure of $(\pm^m|\pm^m)$ hold \cite[Lemma 2.2]{S4};
for example, $Q(\lambda|\pm^m)(\mu|\pm^n)=Q(\lambda\mu|\pm^{m+n})$, etc.

Second, a \textit{$C_a$-labeled} DAG $Q$ is called \textit{slim} if $(b,d)\in Q\Rightarrow (b,d+2)\in Q$ \cite[Definition 2.3]{S4}.
For any slim $Q$, an easy computation gives the inclusion
$\operatorname{Node}Q(\lambda|\pm^m)\hookrightarrow \operatorname{Node}Q(\pm^m|\bar\lambda)_{[-2m]}$ (where $Q_{[d]}$ is $Q$ shifted by depth $d$).
Using this we obtain the new slim DAG $Q[\lambda|\pm^m)$ called (left) \textit{$\lambda$-fragment} of $\mathcal{Q}(\pm^m|\pm^m)$ defined by
\begin{align*}
&\operatorname{Node}Q[\lambda|\pm^m)
:=
\operatorname{Node}Q(\lambda|\pm^m)
\hookrightarrow
\operatorname{Node}Q(\pm^m|\bar\lambda)_{[-2m]},
\\
&\operatorname{Edge}Q[\lambda|\pm^m)
:=
\operatorname{Edge}Q(\lambda|\pm^m)
\cup 
\operatorname{Edge}Q(\pm^m|\bar\lambda)_{[-2m]}|_{\operatorname{Node}Q[\lambda|\pm^m)^2}
\end{align*}
\cite[Definition 2.4]{S4}.
By the (left) \textit{fragmented} $Q(\pm^m|\pm^m)$ we mean $Q[\pm^m|\pm^m):=\sqcup_{\lambda\in\Z_2^m}Q[\lambda|\pm^m)$.
As in the last paragraph, denote $Q[D|E)$ the induced subgraph for  $D,E\subseteq\Z_2^m$, and the composition rules $Q[\lambda|\pm^m)[\mu|\pm^n)=Q[\lambda\mu|\pm^{m+n})$ etc. also hold.
The (right) $\bar\lambda$-fragment $Q(\pm^m|\bar\lambda]$ etc. are defined similarly.
For the simplest slim DAG, i.e. the $C_0$-labeled DAG $[\ |\ )=(\ |\ ]:=2\Z_{\geq 0}$, we abbreviate $[\ |\ )[\lambda|\pm^m)$ to $[\lambda|\pm^m)$, etc.

In the remainder of this subsection we consider in an abelian category $\mathcal{C}$ and identify DAGs with ALDs (and, if no confusion can arise, with the objects or the elements in $[\mathcal{C}]$).
For objects represented by $Q[\lambda_1\pm^m|\pm^{m+1})$ and $Q[\pm^m|\pm^m)$, the two operations $H^0(-)$ and $\slash{\sim}$ are defined by
\begin{align}
\begin{aligned}
\label{two operations}
&H^0(Q[\lambda_1\pm^m|\pm^{m+1}))
:=
Q[\lambda_1\pm^m|+\pm^{m})
=
Q(\lambda_1|+)[\pm^m|\pm^m),
\\ 
&Q[\pm^m|\pm^m)\slash{\sim}
:=
Q[-\pm^{m-1}|\pm^m)
=
Q[-|\pm)[\pm^{m-1}|\pm^{m-1}),
\end{aligned}
\end{align}
respectively.
Denote $\tilde{H}^0(-)$ the composition $\slash{\sim}\circ H^0(-)$.
Then the recursive expression in the third term of 
\begin{align}
\label{nesting for star}
[\lambda_1|+)^\ast[+^{m-1}|-^{m-1})
\simeq
[\lambda_1|+)^\ast[-^{m-1}|+^{m-1})
=
H^0(\tilde{H}^0(\dots(\tilde{H}^0([\lambda_1|\pm)^\ast[\pm^{m-1}|\pm^{m-1}))\dots)
\end{align}
``transports" the first to second.\footnote{In the following, we often consider the cases where only the first parameter is reversed as \eqref{nesting for star}. See \cite[Remark 3.4]{S4} for the reason.}
The nested use of the relations between both sides of \eqref{two operations} in $[\mathcal{C}]$ shows that
\begin{align}
\label{combinatorial formula for star}
\eqref{nesting for star}
=
\sum_{\nu\in\Z_2^m}
(-1)^{|\nu|}
\sum_{n\geq 0}
(\lambda_1\ast\nu_1)\overline{\nu_{2\dots m}}\ast{}[+|\pm)^\ast[+^{m-1}|\pm^{m-1})_{[2n+m+|\nu|-|\overline{\lambda_1}|]}
\binom{n+m-1}{n}
\end{align}
in $[\mathcal{C}]$, where $\lambda\ast Q'$ is $Q'$ shifted by bits $\lambda$ \cite[Section 2.3]{S4}.
At this point, we have obtained a ``combinatorial recursive approximation" of \eqref{RHS of M1, p.10} for the case of star graphs, as detailed and generalized in later subsections.

As mentioned above, in an abelian category there exist degrees of freedom of ALDs due to non-split extensions.
Put another way, we may take a ``linear combination" of nodes of the same color in a (non-slim) DAG and deform its network structure of the arrows.
By taking this advantage, we may consider the ``\textit{defragmentation}"
\begin{align}
\label{defragmentation}
\int_I[\pm^m|\pm^m):=(\pm^I|\pm^I)[\pm^{m-I}|\pm^{m-I})
\end{align}
etc., where $I\subseteq\{1,\dots,m\}$ and $m-I$ the complement (recall the picture of Riemann sum).
We also assume the ``Fubini's theorem" for any $I,J\subseteq\Z_2^m$, $I\cap J=\phi$.
Let us call objects that can be represented in the form $[\lambda|\pm^m)$ etc. up to Cartesian products with some $\mu\ast(\pm^n|\pm^n)_{[d]}$, \textit{type $m$}.
Then the non-split extensions \eqref{defragmentation} enables us to reduce type $m$ objects to $m-|I|$.
In particular, the recursive computation in the last paragraph can be viewed as round-trips between two types: $0$ and $1$ (cf. \cite[Remark 1.2, Section 1.2]{S4}).
And it is the very these two types that constituted the Felder complex and Feigin--Tipunin construction in the case where $\mathcal{C}$ is the module category of Virasoro VOA at level $p\in\Z_{\geq 2}$ corresponding to $3$-leg star graphs.

Actually, the Felder complex and FT construction can be developed in a VOA-independent setting \cite{LS}.
Set
\begin{align*}
[\lambda_1|\pm]
:=
\cdots
\oplus[\lambda_1|\pm)_{[2]}
\oplus[\lambda_1|\pm)
\oplus(\pm|\bar\lambda_1]
\oplus(\pm|\bar\lambda_1]_{[2]}
\cdots,
\ \ 
[\lambda_1|\pm]^h
:=
\begin{cases}
[\lambda_1|\pm)_{[h-|\lambda_1|]}
&h\in |\lambda_1|+2\Z_{\geq 0},
\\
(\pm|\bar\lambda_1]_{[|\lambda_1|-h-2]}
&h\in |\lambda_1|-2\Z_{> 0}.
\end{cases}
\end{align*}
The parameter $h$ is the Cartan weight of $[\lambda_1|\pm]$.
For ease of explanation, first we consider $[\lambda_1|\pm]^\ast$.
Then for the subobject
$[\lambda_1|+]
=
\cdots\oplus[\lambda_1|+)
\oplus(-|\bar\lambda_1]\oplus\cdots$, 
we have the short exact sequence of $H$-modules (Felder complex)
\begin{align}
\label{Felder complex}
0
\rightarrow[\lambda_1|+]^\ast
\rightarrow[\lambda_1|\pm]^\ast
\rightarrow[\bar\lambda_1|+]^\ast(-1)
\rightarrow 0,
\end{align}
where the $(-1)$ means the shift of Cartan weight by $-1$.
In addition, let us assume that $[\lambda_1|\pm]^\ast$ has the negative Borel $B$-action and \eqref{Felder complex} is compatible with the actions. 
Then $[\lambda_1|+]^\ast$ is the maximal $\SL$-submodule of $[\lambda_1|\pm]^\ast$, therefore under the above assumption we may call type $0$ and $1$ objects $\SL$- and $B$-type, respectively.
Furthermore, applying the sheaf cohomology $H^\bullet(\SL\times_B-)$ to \eqref{Felder complex} shows that
\begin{align}
\label{FT construction}
H^0(\SL\times_B[\lambda_1|\pm]^\ast)\simeq[\lambda_1|+]^\ast,
\ \ 
H^1(\SL\times_B[\lambda_1|\pm]^\ast)\simeq 0,
\end{align}
that is the FT construction and the vanishing theorem.
Similar discussions hold for $[\lambda_1|\pm]$ and positive Borel actions etc \cite[Section 3.3]{S4}.
The $H^0(-)$ in \eqref{two operations} and its relation in $[\mathcal{C}]$ are  $H^0(\SL\times_B-)$ and the Weyl-type character formula follows from \eqref{FT construction}, respectively. 
Actually, the discussion of this paragraph can be generalized to any gauge group $G$, and the geometric representation theory of FT construction shows strong similarities with the modular representation theory of $G$ (\cite{LS}; e.g., \eqref{FT construction} is the ``Borel--Weil--Bott theorem").
Also, for the higher rank version of this subsection up to \eqref{Felder complex}, the $B$- or $G$-actions are decorative, as is the VOA-structure; it suffices to replace $\Z_2$ by Weyl/Coxeter groups (cf. \cite[Definition 1.6]{LS}).
It is worth pointing out here that our commitment to abstraction/minimalism is needed not only for the logical transparency, but also to increase the success rate of our aim: to construct a mathematical theory robust enough to hold for all trees, even though we know only the two easiest cases, decorative/contingent elements should be thoroughly separated/eliminated.

Let us summarize this subsection.
For each $m\in\Z_{\geq 0}$, we may consider an abelian category $\mathcal{C}_m$ that consists of the $C_m$-colored type $i\leq m$ objects and is closed under \eqref{two operations} and \eqref{defragmentation} with some compatibilities (such as above ``Fubini's theorem", reversal of the first parameter of type $m$ objects, etc).
As $m$ gets larger, $\mathcal{C}_m$ becomes more complicated, but the semisimple objects \eqref{nesting for star} have the combinatorial relation \eqref{combinatorial formula for star} in $[\mathcal{C}_m]$.
Actually using \eqref{two operations} and \eqref{defragmentation}, this recursion can be considered virtually within $\mathcal{C}_1$, where we can develop the theory of FT construction under the assumption of the $B$-actions, so that \eqref{nesting for star}, \eqref{combinatorial formula for star} are the nested FT constructions/Weyl-type character formula.
In this sense, the study of $\mathcal{C}_m$ is reducible to $\mathcal{C}_1$ or the theory of nested FT, and can be viewed as an abelian categorification of $(m+2)$-leg star graphs or a blueprint for the module category (more precisely, a Weyl group orbit) of the ``Virasoro VOAs". 
If we regard \eqref{nesting for star}/its $H$-extensions as simple objects, then $\mathcal{C}_m$ is for the ``singlet/triplet VOAs" (for the dictionary between ALDs/DAGs and VOAs, see \cite[Section 3.4]{S4}).
It would make sense to call $\mathcal{C}_m$ \textit{hypercubic} or \textit{$m$-cubic}, because the category itself has the nested structure derived from  $(\pm^m|\pm^m)$.
In short, it is an abelian categorical marriage of the basic objects $(\pm^m|\pm^m)$ and $[\ |\ )$.

\section{From star to tree}
\subsection{Argumentation}
\label{argumentation}
Before going into details, let us briefly explain how the discussion should be developed.
Our aim is to provide a reasonable mathematical model for corresponding log VOAs to 3-manifolds.
Therefore we should not only give one such model and deduce $\hat{Z}$, but also show the thought process to get there; this is the reason for the style of this paper.
We first have enough reasons to develop Section \ref{quick review of S4}/nested FT:
as mentioned in Section \ref{introduction}, it not only provides a blueprint for the above correspondence, but also greatly extends the clear perspective/reductive approach that is already successful;
our minimalism maximizes the robustness and scope of the theory;
the observation that the most economical way to extend a procedure is recursion/nesting justifies our model in terms of Occam's razor, etc.
Recalling that rooted trees describe the recursive processes, 
extending Section \ref{quick review of S4} to such general trees seems one of the best options that can be taken at this time.
It is worth noting here that Section \ref{quick review of S4} is independent from $\hat{Z}$.
As already shown in \cite{CCFGH,CFGHP}, the relation between VOA and $\hat{Z}$ is not simple: $\hat{Z}$ are given by \textit{virtual generalized characters} in general, i.e. linear combinations of some $q$-series with the form \eqref{combinatorial formula for star} but with more generalized parameters than VOA-modules allow.
In other words, it is risky to develop Section \ref{quick review of S4}, which is a blueprint for the module categories of the log VOAs, based on  observations of $\hat{Z}$. 
Therefore we first extend Section \ref{quick review of S4} to general trees without reference to $\hat{Z}$, and then by formally substituting the $q$-series $\Theta_{2d}(q)$ in our model, we reconstruct \eqref{RHS of M1, p.10} as a ``virtual generalized character".

\subsubsection{}
\label{argumentanion I}
First, we describe the basic situation that for a given rooted tree $T=(V,E)$, Section \ref{quick review of S4} can be applied at each vertex
according to the tree order (i.e. the processing at a vertex $v$ begins only after all its descendants have been processed), where $v$ is processed means that the recursion (or nested FT) is completed in the sense of \eqref{nesting for star}; we call it \textit{processed/singlet at $v$}.
On the other hand, if no such processing has ever occurred at $v$, we call it \textit{unprocessed/projective at $v$}.
Without loss of generality, we may assume that only nodes are processed.
In fact, by growing leaves so that the degree of each vertex becomes at least $3$, we can change the vertex version to the node version.
Therefore we adopt the node version for later convenience.

It is essentially the semisimplicity of type $0$ objects that makes the recursive process possible.
For example, the state in which all the descendants of $v$ have been processed should seamlessly be considered as the projective object at $v$, allowing the successive applications of Section \ref{quick review of S4} from the descendants to $v$.
What makes it possible is exactly the fact that both singlet (at the descendants) and projective (at $v$) are type $0$.
More generally, if Section \ref{quick review of S4} is applicable at $v$, then the states at the others are type $0$, because the latter information should be embedded into each singlet object at $v$, as in \eqref{defragmentation}, not to affect the hypercubic structure at $v$.
In other words, the existence of multiple nodes means the degrees of failure of the global  hypercubic structure among them. 

Let us express the above considerations more mathematically.
In the remainder of this paper, we consider in an abelian category $\mathcal{C}_T$ associated with $T$.
A number of the recursion $m_v\in\Z_{\geq 1}$ is attached to each $v\in V_{\geq 3}$.
The objects we will consider are labeled by some $(w,X_w)_{w\in V_{\geq 3}}$, where $X_w$ is a type $i\leq m_w$ object.
First we consider the case where for a fixed $v\in V_{\geq 3}$ and a subtree containing all the ascendants of $v$, $X_w$ is a fragment ($v=w$), projective (in the subtree), singlet (otherwise), respectively.
We call such an object a \textit{fragment at $v$}.
For our purpose, we assume that the subcategory consisting of the fragments at $v$ with a fixed state at $w\not=v$ is $m_v$-cubic, and \eqref{two operations}, \eqref{defragmentation} are compatible with the labeling (i.e. objects that have the same label as a result of the operations are isomorphic in $\mathcal{C}_T$). 
The latter assumption justifies, for example, the ``successive applications of Section \ref{quick review of S4}" in the last paragraph; that is, giving the projective object at $v$ by \eqref{defragmentation} (while the descendants are singlet) and giving a singlet object at the descendants of $v$ by \eqref{nesting for star} (while $v$ is projective) produce the same result.
Second we may construct the case where each $X_w$ is a fragment at $w$ as follows.
Let us consider a fragment at some $v$.
This is a subquotient of the projective object at $v$, and the extension of the singlet objects at $v$ described by the ALD structure.
On the other hand, due to the ``successive applications of Section \ref{quick review of S4}" above, these singlet objects at $v$ are the projective objects at an unprocessed sibling of $v$, or if all the siblings have already been processed, at the parent of $v$.
The same as above holds for these nodes, and thus the desired object can be constructed inductively: 
we call such an object a \textit{Fock object of $T$}.\footnote{Unless we fix the order of processing for $V$, we may have to assume that the above inductive constructions are order independent. However, no matter in which order they are constructed, they produce the same result at the level of $[\mathcal{C}_T]$.}

Using these objects, let us generalize \eqref{nesting for star}, \eqref{combinatorial formula for star} to $T$.
For an abelian category $\mathcal{C}$ and $a,b\in[\mathcal{C}]$, denote $a\sim b$ if $[a:m]=0\Leftrightarrow [b:m]=0$ for any simple object $m$, where $[x:m]\in\Z_{\geq 0}$ is the multiplicity of $m$ in $x\in[\mathcal{C}]$.
The \eqref{combinatorial formula for star} was computed in \cite[Section 2.3]{S4} as follows.
The nested use of the relations between both sides of \eqref{two operations} in $[\mathcal{C}]$ first produce the ``projective version" of \eqref{combinatorial formula for star}: the $=$ and fragments $\mu\ast[+|\pm)^\ast[\pm^{m-1}|\pm^{m-1})_{[d]}$ in RHS are replaced by $\sim$ and the corresponding projective object $\mu\ast[\pm|\pm)^\ast[\pm^{m-1}|\pm^{m-1})_{[d]}$, respectively.
However, by some binomial algebraic identity,
it can be replaced by \eqref{combinatorial formula for star}.
It can be extended as follows.
We start from an extreme case: each $X_w$ is singlet at $w$, called a \textit{singlet object of $T$}.
Let us consider the hypercubic structure at the root $\mathrm{rt}(T)$, and then the singlet object is represented by the projective objects at $\mathrm{rt}(T)$ in the sense of the above ``projective version" of \eqref{combinatorial formula for star}.
By repeating the same computation from $\mathrm{rt}(T)$ to leaves, we can finally represent the singlet object of $T$ by the extreme case on the other side: each $X_w$ is projective at $w$, called a \textit{projective object of $T$}, 
up to $\sim$ (if the singlet object of $T$ is considered simple, it is intended to be the projective cover. Other cases, such as ``Virasoro objects", are similarly defined as in \cite[Section 3.4]{S4} and are therefore omitted).
The repeated use of the above binomial algebraic identity allows us to replace the $\sim$ and the projective objects of $T$ by $=$ and the corresponding Fock objects of $T$, respectively.
Needless to say, we may consider the $H$-extensions of the Fock objects of $T$ and can develop the theory of nested FT constructions under the assumption of the $B$-actions.
In any case, with the abbreviated notations $m:=(m_v)_{v\in V_{\geq 3}}$, $[D|E):=([D_v|E_v))_{v\in V_{\geq 3}}$ etc., we have 
\begin{align}
\label{combinatorial formula for tree}
[-|+)^\ast[+^{m-1}|-^{m-1})
=
\sum_{\nu\in\Z_2^m}
(-1)^{|\nu|}
\sum_{n\in\Z_{\geq 0}^{V_{\geq 3}}}
\bar\nu\ast
[+|\pm)^\ast[+^{m-1}|\pm^{m-1})_{[2n+m+|\nu|]}
\prod_{v\in V_{\geq 3}}\binom{n_v+m_v-1}{n_v}
\end{align}
in $[\mathcal{C}_T]$,
where we take $\lambda_1=-$ for simplicity.
This is the tree extension of \eqref{combinatorial formula for star}.

\subsubsection{}
\label{argumentanion II}
The above discussion describes the recursive process when each node has the hypercubic structure, but does not explain why/where the structures arise.
In the vertex version, the most economical (and tautological) solution to this problem is to assume that the processed children of each $v$ are regarded as exactly the parameters of the hypercubic structure at $v$ (where the leaves are considered to be already processed unconditionally).
Since we are considering the node version, this assumption should be slightly modified.
Noting the difference between the root $\mathrm{rt}(T)$ and the other nodes, the most natural modification is the following: 
the recursion number $m_v$ at each node $v$ is $\deg(v)-2$, which consists of one ``$3$-parameter" $\Delta(v)$ that considers some three adjacent vertices of $v$ as one parameter, and the other $\deg(v)-3$ adjacent vertices, called ``$1$-parameters".
Since the parent of $v$ cannot be used as a $1$-parameter of $v$, $\Delta(v)$ consists of the parent and some two children of $v$; we call them \textit{top} and \textit{bottom} of $\Delta(v)$, respectively (in the case $v=\mathrm{rt}(T)$, we call the three children in $\Delta(v)$ the bottom for convenience).
We can also see that since there is only one reversal parameter (negative $B$-action in FT) at each $v$, the $\Delta(v)$ should correspond to it.
Finally, it is natural, if not necessary, to assume that $T$ is a centered tree (because we have no reason to defer the processing of nodes with smaller heights).

Let us take a closer look at the $\{\Delta(v)\}_{v\in V_{\geq 3}}$.
In the following, for the convenience of explanation, we assume that $V_2=\phi$ without loss of generality (we can restore the general case with minor modifications, and in fact $V_2$ plays almost no role in our discussion).
For each node $v\not=\mathrm{rt}(T)$, the top of $\Delta(v)$ has three possibilities: 
$\mathrm{rt}(T)$, 
a $1$-parameter of the grandparent of $v$,
a bottom of the $3$-parameter of the grandparent of $v$.
If we define $\Delta(v)$ and $\Delta(w)$ are \textit{connected} when they share a non-root vertex (this exception of $\mathrm{rt}(T)$ is for later convenience), then $\{\Delta(v)\}_{v\in V_{\geq 3}}$ is decomposed into connected components. 
We denote $\Delta(C)$ the connected component consisting of $\{\Delta(v)\}_{v\in C}$ for some $C\subseteq V_{\geq 3}$, and call the highest node and lowest nodes in $\Delta(C)$ the top and bottom of $\Delta(C)$, respectively.
A connected component has three possibilities: 
the top is $\mathrm{rt}(T)$, 
the top is a $1$-parameter,
it contains $\Delta(\mathrm{rt}(T))$. 
Note that only the last one is unique in general (and only it does not have the top).

With the above additional data on $T$, we can understand the recursive process in \ref{argumentanion I} with higher resolution.
We call a parameter $x$ of some $v$ is \textit{evaluated} if the operation $H^0(-)\circ\slash{\sim}$ is applied to $x$ at $v$. 
In particular, ``$v$ is processed" is equivalent to ``all parameters of $v$ are evaluated".
Using the above structures of connected components, the evaluation can be extended to that of vertices as follows.
If $x\in V$ is a $1$-parameter, the above definition is adopted, otherwise we say that $x$ is evaluated if for the connected component $\Delta(C)$ containing $x$, its top (if there exists) and all $3$-parameters composing $\Delta(C)$ are evaluated. 
We assume that $\mathrm{rt}(T)$ is evaluated unconditionally, because it is the only vertex that is not evaluated in the above definition.
Then the recursive process in \ref{argumentanion I} from a projective object of $T$ to the singlet object of $T$ is equivalent to evaluating all $x\in V$.
(Conversely, as in the last paragraph in \ref{argumentanion I}, when we proceed from a singlet object of $T$ to the projective/Fock objects of $T$, the colors of the latter are determined gradually from $\mathrm{rt}(T)$ to leaves. 
By abuse of notation, we use the same term ``evaluate" for this determination process in the reverse direction.)
Note that a $3$-parameter can be thought of as a ``lazy evaluator" of the vertices composing it, because the final evaluation of them is deferred until the evaluation of the top or  $\Delta(\mathrm{rt}(T))$ of the connected component containing them.
We call a set of $1$- and $3$-parameters given a detail of the lazy evaluation system a \textit{parameter structure on $T$}.

Let us introduce a special case of the parameter structure on $T$.
We consider the above reverse direction: from $\mathrm{rt}(T)$ to leaves.
Suppose that each $x\in V$ is assigned a parameter $e_x\in\{+,-,\pm\}$, where the first two mean the fixed parameters, but the last one means the indeterminate.
The evaluations of $1$- and $3$-parameters determine the parameter of the vertex and the product of the parameters of the  three vertices, respectively.
Following the process in \ref{argumentanion I}, for starting from a singlet object of $T$ with all $e_x$ are $\pm$ and finally arriving at the projective/Fock objects of T with all $e_x$ are $+$ or $-$, the following additional \textit{initial condition} is required: 
the parameter $e_{\mathrm{rt}(T)}$ is fixed from the beginning, and for any $v\in V_{\geq 3}$ the product of the parameters of the bottom vertices of $\Delta(v)$ determines that of each of them.
Once such an initial condition is fixed, for each $v\in V_{\geq 3}$ and the parameter $\nu=e_{\Delta(v)}e_{v_4}\cdots e_{v_{\deg(v)}}$, $\Delta(v):=\{v_1,v_2,v_3\}$,
we may rewrite \eqref{combinatorial formula for tree} as $(-1)^{|\nu|}=\prod_{i=1}^{\deg(v)}e_{v_i}$, etc. 
By taking the sum of \eqref{combinatorial formula for tree} with respect to all the initial conditions, we obtain
\begin{align}
\sum_{e\in\{\pm 1\}^{V_{\geq 3}}}
(\prod_{v\in V_{\geq 3}}e_v^{\deg(v)-|\bar{v}|})
\sum_{\epsilon\in\{\pm 1\}^{V_1}}
(\prod_{i\in\bar{v}}\epsilon_i)
\sum_{n\in\Z_{\geq 0}^{V_{\geq 3}}}
F(n+\tfrac{\deg(v)-2}{2},e,\epsilon)
\prod_{v\in V_{\geq 3}}\binom{n_v+\deg(v)-3}{n_v},
\label{the categorical version of Zhat}
\end{align}
where
$F(n+\tfrac{\deg(v)-2}{2},e,\epsilon)
:=\overline{e_{\Delta(v)}e_{v_4}\cdots e_{v_{\deg(v)}}}\ast[+|\pm)^\ast[+^{\deg(v)-3}|\pm^{\deg(v)-3})_{[2n+\deg(v)-2+|e_{\Delta(v)}e_{v_4}\cdots e_{v_{\deg(v)}}|]}$.
The \eqref{RHS of M1, p.10} is reconstructed by formally substituting $\Theta_{2d}(q)$ in $F$.
However, as noted above, such a substitution may not compatible with our abelian categorical structures  (hence it may be a ``virtual  generalized character").

\subsection{Glimpse of the VOA side}
\label{glimpse of the VOA side}
The above formal substitution is rephrased as a formal interpretation of the Fock object (or the parameter structure on $T$) in terms of characters/conformal weights.
Here we discuss how some of VOA-specific data may be interpreted in our theory.
Recalling that our aim is to study the VOAs from the above abelian categorical perspective, it might be useful to consider not only the VOA-independent sides, but also such VOA data even at a rough level.
Since our theory entirely depends on the tree $T$ and structures on it, we need to give the ``tree expression" of the VOA data.
In the two cases known to date, $3$- or $4$-leg stars \cite{CFGHP}, such a tree expression is given only for the special parameters recovering \eqref{RHS of M1, p.10}. 
However, as it is the special parameter structure in the fourth paragraph of \ref{argumentanion II} (and some additional data), the general parameter structure ($+$ additional data) is regarded as the tree expression of the VOA data for general parameters.
Therefore, for general tree cases, we may consider the natural generalizations of these tree expressions to be the ``VOA data".
This subsection is much more speculative and vague than the discussion up to Section \ref{argumentation}, but may help us to get a rough idea of how the two known cases above extend to our hypothetical log VOAs.

We first recall the case of $4$-leg stars (the $3$-leg case is obtained by replacing $p$ below with $1$). 
The corresponding log VOA is the $(p,p')$-singlet Virasoro VOA for some $p,p'\in\Z_{\geq 2}$, 
$p\perp p'$.
This is an infinite extension of the Virasoro VOA at level $p'/p$, and is embedded into the  Heisenberg VOA (and thus, the lattice VOA $V_{\sqrt{2pp'}\Z}$).
It is worth noting that the ``type of level" determines the Virasoro module structures of the Fock spaces: the cases $p'/p=1$, $p'/p\in\Z_{\geq 2}$, otherwise correspond to type $0,1,2$ objects, respectively.
In that sense, our theory gives a vast expansion/abstraction of the ``type of level/object" of Virasoro VOAs (in general, W-algebras).
The $(p,p')$-singlet has a family of simple modules $W_{s,s'}^\pm$ parametrized by $1\leq s\leq p$, $1\leq s'\leq p'$.
A fixed $(s,s')$ is regarded as a representative of the $\Z_2$-orbit with respect to the $\Z_2$-actions
$W_{s,s'}^\pm\leftrightarrow W_{p-s,s'}^\mp$ and  $W_{s,s'}^\pm\leftrightarrow W_{s,p'-s'}^\mp$.
The module category is decomposed into these $\Z_2$-orbits, and our $\mathcal{C}_T$ (now, $\mathcal{C}_2$) is intended to be one of them.

We first consider the ``levels" and ``representatives of simple modules" in our theory following the first paragraph.
The \eqref{RHS of M1, p.10} and $W_{s,s'}^\pm$ are related as follows:
for some $p_1,\dots,p_4\in\Z_{\geq 2}$ associated with each leg of the plumbed graph, the \eqref{RHS of M1, p.10} are given by the virtual generalized characters of $W_{s,s'}^\pm$ with $p=p_4$, $p'=p_1p_2p_3$, $s=1$, and $s'$ with a form $s'_{\lambda_1,\lambda_2,\lambda_3}$ ($\lambda_i\in\Z_2$), where which $p_i$ is assigned as $p=p_4$ is arbitrary \cite[Section 4.2]{CFGHP}.
This specialization of levels/representatives resulting from too many parameters in the plumbed graph fits into the special case in the fourth paragraph of \ref{argumentanion II}:
the $1$- and $3$-parameters are the leaves corresponding to $p_4$ and $p_1,p_2,p_3$, and for the initial conditions $(\lambda_i)_{1\leq i\leq 3}$, the representatives are defined by $(s,s')=(1,s'_{\lambda_1,\lambda_2,\lambda_3})$.
Considering inversely, the general parameter structure on $T$ is the ``tree expression" of the level/representative for general $p,p',s,s'$;
that is, we simply associate $(p,s)$ and $(p',s')$ with the $1$- and $3$-parameters, respectively.
Therefore, we may expect in general that the parameter structure and the partial data $\{\Delta(v)\}_{v\in V_{\geq 3}}$ (it automatically determines $1$-parameters too) are regarded as the ``type of representative" and the ``type of level", respectively, and some numerical input $(\vec{p},\vec{s})$ (resp. $\vec{p}$) to them  reconstructs the actual ``representative"  (resp. ``level"). 
In the special case, $\vec{p},\vec{s}$ are determined by $(w_v)_{v\in V}$ and the initial condition.

We second consider the  ``conformal weights" in our theory.
For this, we have to associate Fock objects with a positive definite lattice.
As we take the recursive process at each node separately, we may use a rank $|V_{\geq 3}|$ lattice.
Let us restrict ourselves to the special parameter structure on general $T$.
Recalling the $3$- or $4$-leg cases, each $x\in V_1$ seems to give the fractional part $e_x/2w_x$ in the $v$-th dimension such that $x\in\bar{v}$ (the fractional part gives the above representative $(1,s_{\lambda_1,\lambda_2,\lambda_3})$).
On the other hand, the simplest interpretation of $e_v$ and the depth at $v\in V_{\geq 3}$ is that they are the orthant (i.e. positive or negative direction/half-line) and the absolute value in the $v$-th dimension, respectively (and thus, the recursive process from $\mathrm{rt}(T)$ to leaves gradually determines the coordinates of the lattice).
Note that the above formal substitution is to some extent reasonable in the sense that by comparing \eqref{RHS of M1, p.10} with the special parameter structure, the same parameter interpretation is obtained.
Therefore, the ``conformal weights" are considered to be reconstructed from the above ``representatives" and the parameter interpretation (the quadratic form is determined by the partial data: plumbed graph).
Similarly, for general parameter structures on $T$, the ``representatives" and generalized parameter interpretations (i.e. the coordinate in the $v$-th dimension is determined when $v$ is evaluated) are expected to be the desired data.

As in the case of $(p,p')$-singlet, our general tree cases are also expected to be the extensions of the ``Virasoro VOAs" (in general, ``W-algebras") with the above greatly extended levels and parameters. 
Recalling that in the case of $3$-leg stars, the FT construction combined with the above three VOA data proved the duality/simplicity of the Arakawa--Frenkel modules \cite{ArF} of the corresponding W-algebras \cite{CNS,LS,S2}, our ``W-algebras" may provide vast generalizations for which such an approach is valid (cf. \cite[Remark 1.2]{S4}).

\subsection{Discussions and future directions}
\label{conclusion}
We first conclude the above discussion.
The type of level (or that of Fock objects) is regarded as the ``\textit{distance from semisimplicity}" of our abelian category. 
Therefore, \cite{Zhu} shows that our recursion theory is very natural and straightforward for the construction of the (abelian category of modules of the hypothetical) log VOAs. 
The binary nature of our distance leads to the hypercubic structures, and thus the character formula \eqref{combinatorial formula for tree} for simple objects in terms of Fock objects is obtained by using recursion between types 0 and 1. 
Moreover, when each step has an additional $B$-action, it is obtained in the geometric representation theory of nested FT  constructions,
and on the other hand, by ``averaging" the characters with respect to the initial conditions, the corresponding $\hat{Z}$-invariant is reconstructed \eqref{the categorical version of Zhat} up to $\Theta_{2d}(q)$. 
Comparing the discussion with \eqref{Zhat-invariant}, we can see that our theory is the algebraic counterpart of the contribution from $F_{3d}(x)$ (it also provides a ``physical explanation" for the VOA-independence of our theory: because it essentially belongs to the 3d $\mathcal{N}=2$ theory, rathar than 2d CFT). 
In other words, to investigate the 3d-3d correspondence, it seems that the ``core" is not the log VOA itself but the  combinatorial/Lie algebraic patterns that reflect its \textit{complexity}.

Future directions include the more advanced geometric representation theory of (nested) FT constructions, as well as the implementation of log VOAs/W-algebras. 
However, a much more interesting question is: 
our theory can be extended beyond the negative definite plumbed 3-manifolds?
As far as the author knows, the general explicit formulas for $\hat{Z}$ are currently known only for these cases, but the following two extensions have already been calculated and studied for several examples (e.g. \cite{CCFGH,CFGHP,CDGG,GM}).
The first extension considers the \textit{orientation reversing} of 3-manifolds.
In the $q$-series side, this is expected to be given by ``$q\leftrightarrow 1/q$ duality"  (however, it does not mean simply replacing each $q^n$ by $q^{-n}$). 
The author expect that extending the Kazama--Suzuki coset construction $W_{1/k}(\mathfrak{sl}_{2|1})\simeq \operatorname{Com}(\pi, V_k(\mathfrak{sl}_2)\otimes V_{\Z})$ to the generalized ``levels" in the sense of Section \ref{glimpse of the VOA side} may be a promising direction. 
As mentioned above, our nesting approach is expected to be valid for the representation theory of $V_k(\mathfrak{sl}_2)$, and thus for $W_{1/k}(\mathfrak{sl}_{2|1})$ so is (cf. \cite{CNS}).
The second extension considers \textit{hyperbolic} 3-manifolds. 
At this point, the author does not have any specific ideas on this, but recalling the last paragraph, non-$C_2$-cofiniteness of hyperbolic cases \cite{CCFGH}, and the JSJ decomposition of 3-manifolds, we can pose the following speculative, vague, but thought-provoking question: \textit{Is the log VOA corresponding to a 3-manifold $M$ essentially reducible to a rational VOA by finitely nested FT constructions if and only if $M$ is non-hyperbolic?}

\end{document}